\documentclass[12pt]{amsart}
\usepackage{amsfonts,amsmath,amsthm,amssymb}
\usepackage{latexsym}
\usepackage{graphics}
\newtheorem{theorem}{Theorem}[section]
\newtheorem{corollary}[theorem]{Corollary}

\newtheorem{lemma}[theorem]{Lemma}

\newtheorem{fact}{Fact}

\newtheorem{proposition}[theorem]{Proposition}

\theoremstyle{definition}
\theoremstyle{remark}
\newtheorem{rem}[theorem]{Remark}
\theoremstyle{remark}

\newcommand{\beql}[1]{\begin{equation}\label{#1}}

\newcommand{\eeq}{\end{equation}}

\begin{document}

\title[A two-parameter family of complex Hadamard matrices of order $6$]{A two-parameter family of complex Hadamard matrices of order $6$ induced by hypocycloids}

\author{Ferenc Sz\"oll\H{o}si}
\date{March, 2009.}
\address{Ferenc Sz\"oll\H{o}si: Central European University (CEU), Institute of Mathematics and its Applications, H-1051, N\'ador u. 9., Budapest, Hungary}\email{szoferi@gmail.com}
%\thanks{Supported by Otka.......}

\begin{abstract}
Constructions of Hadamard matrices from smaller blocks is a
well-known technique in the theory of real Hadamard matrices:
tensoring Hadamard matrices and the classical arrays of
Williamson, Ito are all procedures involving smaller order
building blocks. We apply a new block-construction for order $6$
to obtain a previously unknown $2$-dimensional family of \emph{complex} Hadamard matrices. Our results extend the families
$D_6(t)$ and $B_6(\theta)$ found by various authors recently
\cite{BN}, \cite{dita}. As a direct application the existence of a
\emph{$2$-parameter family of  MUB-triplets} of order $6$ is shown.
\end{abstract}
\maketitle {\bf 2000 Mathematics Subject Classification.} Primary
05B20, secondary 46L10.

\section{Introduction}\label{sec:intro}
Constructions of complex Hadamard matrices of small orders was
originally motivated by a question of Enflo, who asked whether for
prime orders an enphased and permuted version of the Fourier
matrix is the only circulant matrix whose column vectors are
bi-unimodular. It was known that this is true for $p=2,3$ but a
subsequent general construction due to Bj\"orck \cite{BJ} (see
also the papers written by Munemasa and Watatani \cite{MW}, and by
de la Harpe and Jones \cite{HJ}) showed that there are
inequivalent examples already for any prime $p\geq 7$. The only
remaining case $p=5$ was settled by Haagerup who has \emph{fully
classified} complex Hadamard matrices up to order $5$, and showed
that Enflo's hypothesis is still true for $p=5$. Haagerup also
pointed out some possibilities for parametrization in composite
dimensions, and introduced an \emph{invariant set} in order to
distinguish inequivalent complex Hadamard matrices from each other
\cite{haagerup}. Currently the smallest order where full
classification is not available is order $6$. Another significant
paper on complex Hadamard matrices was the \emph{cathaloge of
complex Hadamard matrices of small orders} by Tadej and
\.Zyczkowski who, besides introducing another invariant, the
\emph{defect}, listed all known parametric families of complex
Hadamard matrices up to order $16$ \cite{karol}. Most of the
presented matrices could be obtained via Di\c{t}\u{a}'s general
method \cite{dita}, but matrices due to Bj\"orck \cite{BJ},
Nicoara, Petrescu \cite{petrescu} and Tao \cite{tao} have also
been exhibited. Recently the online version of this cathaloge has
significantly been extended by new matrices in at least two
different ways: firstly a new general construction of Butson-type
matrices (i.e.\ matrices built from roots of unities) was
discovered by Matolcsi, R\'effy and Sz\"oll\H{o}si \cite{MM}, who
used a spectral set construction from \cite{KM}, while another
independent construction of Sz\"oll\H{o}si showed how to introduce
parameters to real Hadamard- and real conference matrices to
obtain parametric families of complex Hadamard matrices \cite{sz}.
Secondly, new order $6$ matrices were constructed by Beauchamp and
Nicoara \cite{BN} and by Matolcsi and Sz\"oll\H{o}si \cite{MS}. In
particular all self-adjoint complex Hadamard matrices have been
classified, and a family of symmetric matrices has been
introduced, respectively. On the one hand, constructing complex
Hadamard matrices of order $6$ is interesting of its own as
currently this is the smallest order where full classification of
Hadamard matrices is not available. While recent numerical
evidence suggests that the set of complex Hadamard matrices of
order $6$ forms a $4$-dimensional manifold \cite{skinner}, it
seems that describing all of them through closed analytic
formul\ae\ remains elusive. On the other hand, complex Hadamard
matrices play an important r\^{o}le in the theory of operator
algebras \cite{popa}, and also in quantum information theory
\cite{wer}. In particular, the question whether there exist $d+1$
mutually unbiased bases (MUBs) in $\mathbb{C}^d$ is equivalent to
the existence of certain complex Hadamard matrices, as such bases
can always be taken as a union of the identity operator and a set
of (rescaled) complex Hadamard matrices whose normalized product
is also a (rescaled) Hadamard. For a survey on the MUB problem see
e.g.\ \cite{MUB6}, while for a comprehensive list of applications
we refer the reader to the recent book of Horadam \cite{Ho}.

Our paper is organized as follows: in section $2$ we derive a
\emph{two-parameter family} of complex Hadamard matrices of order
$6$ by considering circulant block-matrices of order $3$. In
section $3$ we discuss some connections between this new family
and some other previously known examples of Hadamard matrices. In
particular, we show that besides some well-known matrices such as
$C_6$ and the members of the generalized Fourier families
$F_6(1,3)$, $F_6^T(1,3)$, the whole affine family $D_6(t)$ and all
self-adjoint Hadamard matrices of order $6$, denoted by
$B_6(\theta)$, belong to our family. All the mentioned matrices
can be found online at \cite{web}. In the last section we recall a
construction of Zauner  \cite{Z} to prove the existence of a {\it
two-parameter family of MUB-triplets} of order $6$. The main
ingredient to his construction is essentially a $2$-circulant
complex Hadamard matrix.
\section{The construction}
The main idea of our method is to consider Hadamard matrices with
a ``highly symmetrical'' block structure. Such restrictions made
on the matrix implies that ``almost all'' orthogonality conditions
immediately hold. We begin our construction with the following
$2p\times 2p$ matrix consisting $p\times p$ blocks of matrices $A,
B$ and their adjoints $A^\ast,B^\ast$ respectively. \beql{}
H=\left[
\begin{array}{cc}
A & B\\
B^\ast & -A^\ast\\
\end{array}
\right]. \eeq In order to $H$ be a complex Hadamard,  one must
exhibit certain unimodular matrices $A,B$ satisfying the following
conditions: \beql{2} AA^\ast+BB^\ast=2pI_p \eeq \beql{4} B^\ast
B+A^\ast A=2pI_p \eeq \beql{3} AB-BA=0, \eeq
where $I_p$ is the identity matrix of order $p$.  Observe that if
we choose $A$ and $B$ to be \emph{circulant} matrices they will
commute, and therefore \eqref{3} will hold identically, while
\eqref{4} will be equivalent to \eqref{2}. Hence, by considering
$p=3$ the building blocks of $H$ can be taken as \beql{} A=\left[
\begin{array}{ccc}
a & b & c\\
c & a & b\\
b & c & a\\
\end{array}
\right],
B=\left[
\begin{array}{ccc}
d & e & f\\
f & d & e\\
e & f & d\\
\end{array}
\right],
\eeq
and so we have $H$ and its dephased form, $X_6$
\beql{8}
H=\left[
\begin{array}{ccc|ccc}
 a & b & c & d & e & f \\
 c & a & b & f & d & e \\
 b & c & a & e & f & d \\
 \hline
 \frac{1}{d} & \frac{1}{f} & \frac{1}{e} & -\frac{1}{a} & -\frac{1}{c} & -\frac{1}{b} \\
 \frac{1}{e} & \frac{1}{d} & \frac{1}{f} & -\frac{1}{b} & -\frac{1}{a} & -\frac{1}{c} \\
 \frac{1}{f} & \frac{1}{e} & \frac{1}{d} & -\frac{1}{c} & -\frac{1}{b} & -\frac{1}{a}
\end{array}
\right],
X_6=\left[
\begin{array}{cccccc}
 1 & 1 & 1 & 1 & 1 & 1 \\
 1 & \frac{a^2}{b c} & \frac{a b}{c^2} & \frac{a f}{c d} & \frac{a d}{c e} & \frac{a e}{c f} \\
 1 & \frac{a c}{b^2} & \frac{a^2}{b c} & \frac{a e}{b d} & \frac{a f}{b e} & \frac{a d}{b f} \\
 1 & \frac{a d}{b f} & \frac{a d}{c e} & -1 & -\frac{a d}{c e} & -\frac{a d}{b f} \\
 1 & \frac{a e}{b d} & \frac{a e}{c f} & -\frac{a e}{b d} & -1 & -\frac{a e}{c f} \\
 1 & \frac{a f}{b e} & \frac{a f}{c d} & -\frac{a f}{c d} & -\frac{a f}{b e} & -1
\end{array}
\right]. \eeq Let us recall that two complex Hadamard matrices,
$H$ and $K$, are called \emph{equivalent}, if there exists $D_1,
D_2$ unitary diagonal and $P, Q$ permutational matrices, such that
$H=D_1PKQD_2$.
\begin{rem}\label{remx}
Clearly, we are free to set $a=d=1$ by natural equivalence. Also, observe that whenever $b=\overline{c}$ we get a self-adjoint matrix by construction. The same holds for the case $e=\overline{f}$ too, as the r\^ole of the blocks $A, B$ are symmetric under the equivalence.
\end{rem}
As we have imposed the circularity conditions on the building
blocks of $H$, \eqref{2} is the only equation to be satisfied. In
other words, it is necessary and sufficient for $H$ to be a
Hadamard matrix to find unimodular complex numbers $a,b,c,d,e,f$,
such that \beql{x}
\frac{a}{b}+\frac{b}{c}+\frac{c}{a}+\frac{d}{e}+\frac{e}{f}+\frac{f}{d}=0
\eeq holds. At the first glance it seems that we have \emph{so
much freedom} to choose $b,c,e,f$ to satisfy \eqref{x}, however,
later we will see that there is a really strong connection between
these seemingly free parameters. Nevertheless we will fully
classify this type of matrices obtaining a new $2$-dimensional
family. Let us denote by $\varphi[x,y] :
\mathbb{T}\times\mathbb{T}\to \mathbb{C}$ the following
fundamental function of ours: \beql{}
\varphi[x,y]:=x+y+\frac{1}{xy}. \eeq Now observe, that we have
$\varphi[\frac{a}{b},\frac{b}{c}]=\frac{a}{b}+\frac{b}{c}+\frac{c}{a}$.
Hence to satisfy \eqref{x} one should look for certain $x,y,u$ and
$v\in\mathbb{T}$, such that for some $\alpha\in
\mathrm{ran}\varphi$ \beql{} \varphi[x,y]=\alpha \eeq and \beql{}
\varphi[u,v]=-\alpha \eeq hold simultaneously. Therefore we should
understand the range of $\varphi$ and characterize the set
$\mathrm{ran}\varphi\cap\mathrm{ran}(-\varphi)$. We recall the
following well-known
\begin{fact}
$\varphi[x,x]=2x+\frac{1}{x^2}$ is a special plane algebraic
curve, a three-sided hypocycloid, called \emph{deltoid}.
\end{fact}
The following is also relatively easily seen.
\begin{fact}
For any fixed $y_0\in\mathbb{T}$, $\varphi[x,y_0]$ is a sliding
line segment with each end on the deltoid and tangent to the
deltoid. Therefore $\varphi[x,y]$ is the union of all such line
segments, i.e.\ the whole interior of the deltoid.
\end{fact}
\begin{figure}
  \caption{The intersection of the two deltoids is the fundamental region $\mathbb{D}$.}
  \begin{center}
   \includegraphics{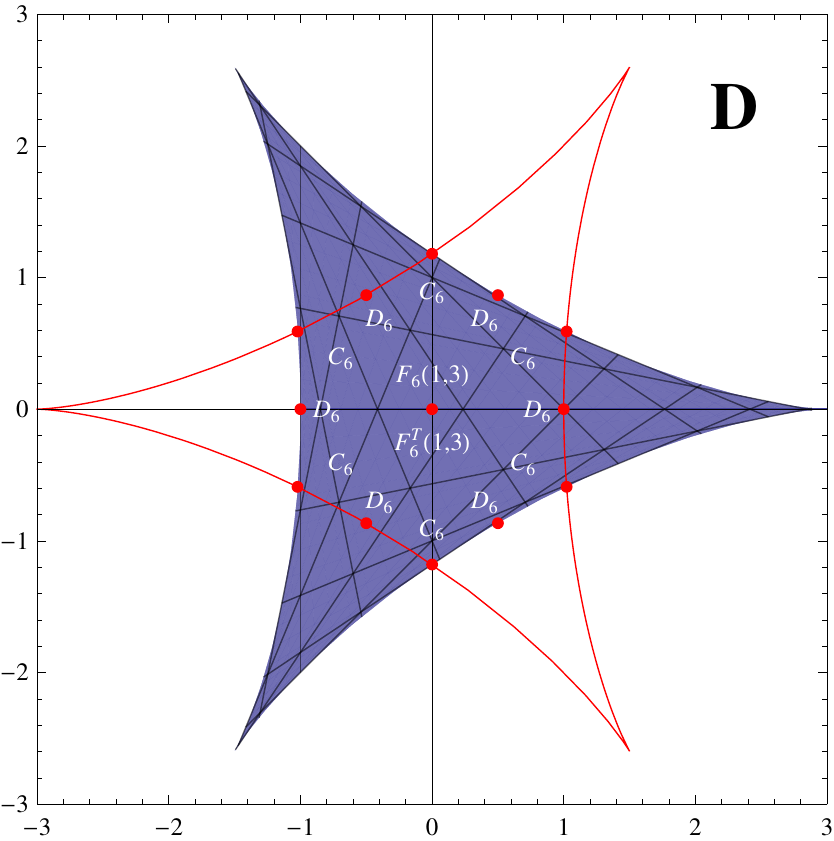}
  \end{center}
\end{figure}
Let us denote the intersection of the two deltoids above by
$\mathbb{D}:=\mathrm{ran}\varphi\cap\mathrm{ran}(-\varphi)$. It is
clear that for any $\alpha\in\mathbb{D}$ one can define a complex
Hadamard matrix in the following way: take any value of
$\varphi^{-1}[\alpha]$, say $x,y$. Then we have $a=1,
b=\overline{x},c=\overline{x}\overline{y}$. Similarly, take
$\varphi^{-1}[-\alpha]$ to obtain the values of $u,v$. We have
$d=1, e=\overline{u},f=\overline{u}\overline{v}$. In particular, we have \beql{hxyuv} X_6(\alpha)\equiv
X_6(x,y,u,v)= \left[
\begin{array}{cccccc}
 1 & 1 & 1 & 1 & 1 & 1 \\
 1 & x^2 y & x y^2 & \frac{x y}{u v} & u x y & v x y \\
 1 & \frac{x}{y} & x^2 y & \frac{x}{u} & \frac{x}{v} & u v x \\
 1 & u v x & u x y & -1 & -u x y & -u v x \\
 1 & \frac{x}{u} & v x y & -\frac{x}{u} & -1 & -v x y \\
 1 & \frac{x}{v} & \frac{x y}{u v} & -\frac{x y}{u v} & -\frac{x}{v} & -1
\end{array}
\right]. \eeq In the rest of this section we describe an algebraic
way of inverting $\varphi$, i.e.\ how we can determine $x,y$ and
$u,v$ from a given $\alpha\in\mathbb{D}$. Considering the equation
$\varphi[x,y]=\alpha$, we have
\beql{122}x+y+\frac{1}{xy}=\alpha.\eeq
After conjugating and using the fact that $x,y,\in\mathbb{T}$
we have
\beql{132}\frac{1}{x}+\frac{1}{y}+xy=\overline{\alpha}\eeq
Instead of solving the system of equations \eqref{122}--\eqref{132} we multiply equation \eqref{122} by $x^2\neq0$ and \eqref{132} by $x\neq0$ and rather consider their sum and difference respectively. In this way the variable $y$ vanishes from the difference and we obtain the following \emph{cubic} equation
for $x$, depending on $\alpha$.
\beql{112}f_\alpha(x):=x^3-\alpha x^2+\overline{\alpha}x-1=0\eeq
It is important to realize that $y$ is a root of \eqref{112} by symmetry as well. Moreover, if $x$ and $y$ are distinct roots of \eqref{112} then \eqref{122} follows. For $\alpha\in\mathrm{int}\mathbb{D}$ the roots of \eqref{112} are
distinct, and let us denote them by $r_1,r_2,r_3$. Our construction guarantees that two of them are unimodular. But as $r_1r_2r_3=1$ we conclude that the third root is unimodular as well. Hence one can choose $x$ as any of $r_1,r_2,r_3$, and choose $y$ as \emph{any
other} root. We therefore have $6$ choices for the ordered
pair $(x,y)$.

Finally, let us substitute $-\alpha$ into \eqref{112}, and denote
the roots by $q_1, q_2, q_3$. The method to determine the values
of $u,v$ is completely analogous to what we have presented for $x$
and $y$.

For $\alpha\in\mathrm{int}\mathbb{D}$ we therefore have $6\times
6=36$ choices for the ordered quadruple $(x,y,u,v)$. However an
easy automatized calculation shows that all of the emerging
matrices $X_6(x,y,u,v)$ are equivalent to one of the two matrices
$X_6(r_1, r_2, q_1, q_2)$ or $X_6^T(r_1, r_2, q_1, q_2)$ (note
that a complex Hadamard matrix is generically not equivalent to
its transpose).\footnote{We thank Ingemar Bengtsson for pointing out
that $X_6(x,y,u,v)$ and $X_6^T(x,y,u,v)$ are generically not
equivalent.} On the boundary of $\mathbb{D}$, however, it is easy
to show that the roots of \eqref{112} are $r, r$ and $\frac{1}{r^2}$ and the two families $X_6(r_1, r_2, q_1, q_2)$ and $X_6^T(r_1, r_2, q_1,
q_2)$ are equivalent, and hence in this case all choices of the quadruple
$(x,y,u,v)$ lead to equivalent matrices.

Finally we note that for every $\alpha\in\mathbb{D}$ $X_6(\alpha)$ is stable under conjugation, that is $X_6(r_1, r_2, q_1, q_2)$ and $\overline{X}_6(r_1, r_2, q_1, q_2)$ are equivalent.

By summarizing the contents of this section, we establish the main
result of the paper.
\begin{theorem}
There exist two previously unknown $2$-parameter non-affine
families of complex Hadamard matrices of order $6$, $X_6(x,y,u,v)$
and $X_6^T(x,y,u,v)$, described by formula \eqref{hxyuv} and its
transposed. For $\alpha\in\mathbb{D}$ the values of $x$ and $y$
are determined as roots of $f_{\alpha}$ in \eqref{112}, while the
values of $u$ and $v$ are determined as roots of $f_{-\alpha}$ in
\eqref{112}, respectively.
\end{theorem}
\begin{proof}
The construction above lead us to a $2$-parameter family of
complex Hadamard matrices as follows. For $\alpha\in\mathbb{D}$
let $r_1(\alpha), r_2(\alpha), r_3(\alpha)$ denote the roots of
equation \eqref{112}, being set as continuous functions of
$\alpha$. For a given $\alpha\in\mathbb{D}$ one can set
$x=r_1(\alpha)$ and $y=r_2(\alpha)$. Similarly, substitute
$-\alpha$ into \eqref{112} and denote the roots as $q_1(\alpha),
q_2(\alpha), q_3(\alpha)$, and set $u=q_1(\alpha)$ and
$v=q_2(\alpha)$. Finally, define $X_6(\alpha)=X_6(x,y,u,v)$ as in
formula \eqref{hxyuv}. We emphasize again that easy permutation
equivalences show that all choices of the roots
$r_{i_1}(\alpha),r_{i_2}(\alpha),q_{i_1}(\alpha),q_{i_2}(\alpha)$
lead to matrices equivalent to $X_6(r_1,r_2,q_1,q_2)$ or
$X_6^T(r_1,r_2,q_1,q_2).$

The main claim of the Theorem is that this family (and its
transposed) has not appeared in the literature so far. To show
this, recall that with the exception of the Fourier families
$F_6(a,b)$ and $F_6^T(a,b)$, all previously known families of
order $6$ contain less than two parameters. Therefore we only need
to exhibit one particular matrix from our family $X_6(\alpha)$
which does not belong to the Fourier families. Such a matrix can
be obtained by choosing e.g.\ $\alpha_0=1$ on the boundary of
$\mathbb{D}$. It is easy to show that in this case all choices of
$(x,y,u,v)$ lead to a Hadamard matrix equivalent to $D_6$, which
is not included in the families $F_6(a,b)$ and $F_6^T(a,b)$.
Therefore, by continuity, in a small neighborhood $\mathcal U$ of
$\alpha_0=1$ the family $X_6(\alpha)$ is disjoint from $F_6(a,b)$
and $F_6^T(a,b)$. Hence, inside this neighbourhood $\mathcal U$
only one-parameter curves can possibly produce already known
complex Hadamard matrices of order 6, while generically
$X_6(\alpha)$ is indeed new.
\end{proof}
This shows that the family $X_6(\alpha)$ is at least {\it locally}
new, around $\alpha_0=1$. We expect that more is true: the family
$X_6(\alpha)$ intersects the Fourier family only at $\alpha=0$.
\section{Connections to previously discovered families}
In this section we analyze how the obtained new family of complex
Hadamard matrices $X_6$ is related to the previously discovered ones,
such as $B_6, D_6, F_6, M_6$ and $S_6$, respectively. In particular, we
prove that both the Beauchamp--Nicoara family of self-adjoint
complex Hadamard matrices and Di\c{t}\u{a}'s one-parameter affine
family is contained in the orbit of $X_6(\alpha)$. Thus our
construction in some sense \emph{unifies and extends} some of the
previously discovered families.

We shall denote the standard basis of $\mathbb{C}^6$ by $e_i,
i=1,2,\hdots,6$, which should be understood as \emph{column
vectors}. Also, for later purposes let us denote by $D[\alpha]$
the discriminant function associated to \eqref{112}, i.e.\ let
$r_1,r_2,r_3$ be the three roots of $f_\alpha$ and define
\beql{d1}
D[\alpha]:=(r_1-r_2)^2(r_2-r_3)^2(r_3-r_1)^2=|\alpha|^4+18|\alpha|^2-8\Re[\alpha^3]-27.
\eeq Clearly, $D[\alpha]\in\mathbb{R}$, and $\alpha\in\mathbb{D}$ if and only if $D[\alpha]\leq0$ and $D[-\alpha]\leq0$. Note also, that on the
boundary of $\mathbb{D}$ we have $D[\alpha]=0$ or $D[-\alpha]=0$.

We begin our investigation with the center of $\mathbb{D}$, i.e.\
we consider the case $\alpha=0$. We have the following
\begin{lemma}
For $\alpha=0$ one choice of $(x,y,u,v)$ in formula \eqref{hxyuv}
leads to a Hadamard matrix equivalent to $F_6(1,3)$.
\end{lemma}
\begin{proof}
Straightforward computation.
\end{proof}

Next we classify the ``extremal'' points of $\mathbb{D}$. It has
six points which are farthest from the center, and another six
which are closest to it. These points will be called ``maximal''-
and ``minimal'' extremal points of $\mathbb{D}$.
\begin{lemma}
a) The six maximal extremal points of $\mathbb{D}$ can be obtained
by choosing \beql{}
\alpha^{max}_k=\sqrt{-9+6\sqrt{3}}\mathbf{e}^{\mathbf{i}\left(\frac{\pi}{6}+k\frac{\pi}{3}\right)},
k=1,2,\hdots, 6 \eeq and lead to matrices equivalent to $C_6$.

b) The six minimal extremal points can be obtained by choosing
\beql{} \alpha^{min}_k=\mathbf{e}^{\mathbf{i}k\frac{\pi}{3}},
k=1,2,\hdots, 6 \eeq and lead to matrices equivalent to $D_6$.
\end{lemma}
\begin{proof}
Straightforward computation.
\end{proof}
Somewhat surprisingly it turns out that the whole family $D_6(t)$
is included in in our family $X_6(\alpha)$. This was actually
first found by Zauner \cite{Z}. We have the following
\begin{proposition}[cf.\ Ex.\ 5.7.\ from \cite{Z}]
Let $D(t)$ be a complex Hadamard matrix of the form \beql{} D(t)=
\left[
\begin{array}{cccccc}
 1 & 1 & 1 & 1 & 1 & 1 \\
 1 & -1 & -\frac{\mathbf{i}}{t^3} & \mathbf{i} & -\mathbf{i} & \frac{\mathbf{i}}{t^3} \\
 1 & -\mathbf{i} t^3 & -1 & -\mathbf{i} & \mathbf{i} t^3 & \mathbf{i} \\
 1 & \mathbf{i} & -\mathbf{i} & -1 & \mathbf{i} & -\mathbf{i} \\
 1 & -\mathbf{i} & \frac{\mathbf{i}}{t^3} & \mathbf{i} & -1 & -\frac{\mathbf{i}}{t^3} \\
 1 & \mathbf{i} t^3 & \mathbf{i} & -\mathbf{i} & -\mathbf{i} t^3 & -1
\end{array}
\right],
\eeq
where $t\in\mathbb{T}$ is an indeterminate. Then $D(t)$ has a $2$-circulant representation.
\end{proposition}
\begin{proof}
Let us define the unitary diagonal  matrices
$D_1=\mathrm{Diag}\left(1,\mathbf{i}t,\mathbf{i}/t,1,t,-1/t\right)$
and
$D_2=\mathrm{Diag}\left(1,\mathbf{i}/t,\mathbf{i}t,1,1/t,-t\right)$.
Then one gets \beql{f32} D_1D(t)D_2= \left[
\begin{array}{ccc|ccc}
 1 & \frac{\mathbf{i}}{t} & \mathbf{i} t & 1 & \frac{1}{t} & -t \\
 \mathbf{i} t & 1 & \frac{\mathbf{i}}{t} & -t & 1 & \frac{1}{t} \\
 \frac{\mathbf{i}}{t} & \mathbf{i} t & 1 & \frac{1}{t} & -t & 1 \\
 \hline
 1 & -\frac{1}{t} & t & -1 & \frac{\mathbf{i}}{t} & \mathbf{i} t \\
 t & 1 & -\frac{1}{t} & \mathbf{i} t & -1 & \frac{\mathbf{i}}{t} \\
 -\frac{1}{t} & t & 1 & \frac{\mathbf{i}}{t} & \mathbf{i} t & -1
\end{array}
\right].
\eeq
\end{proof}
\begin{corollary}
All members of the Di\c{t}\u{a}-family $D_6(t)$ have a $2$-circulant
representation.
\end{corollary}
\begin{proof}
The family $D(t)$ above is trivially permutation equivalent to the
Di\c{t}\u{a}-family $D_6(t^3)$ as listed in \cite{karol}.
\end{proof}

Next we turn our attention to the family of self-adjoint complex
Hadamard matrices.% of order $6$. We have the following
\begin{lemma}\label{salemma}
On the boundary of $\mathbb{D}$ all emerging matrices are self-adjoint.
\end{lemma}
\begin{proof}
Let us suppose that $\alpha\in\partial\mathbb{D}$. Then we have
either $D[\alpha]=0$ or $D[-\alpha]=0$. Suppose that
$D[\alpha]=0$. We have already seen that in this case the roots of \eqref{122} are $r,r$ and $\frac{1}{r^2}$, and as we obtain equivalent matrices we are free to set $x=r$, $y=\frac{1}{r^2}$. Then, as we have $b=\overline{x}$, $c=\overline{x}\overline{y}$, the statement follows from Remark \ref{remx}. The case $D[-\alpha]=0$ is completely analogous.
\end{proof}

It turns out that {\it all} complex self-adjoint Hadamard matrices
of order 6 have a 2-circulant representations.
\begin{proposition}
Let $B$ be a complex Hadamard matrix of the form
\beql{ccc}
B=\left[
\begin{array}{cccccc}
 1 & 1 & 1 & 1 & 1 & 1 \\
 1 & -1 & -\frac{1}{x} & -y & y & \frac{1}{x} \\
 1 & -x & 1 & y & \frac{1}{z} & -\frac{1}{x y z} \\
 1 & -\frac{1}{y} & \frac{1}{y} & -1 & -\frac{1}{x y z} & \frac{1}{x y z} \\
 1 & \frac{1}{y} & z & -x y z & 1 & -\frac{1}{x} \\
 1 & x & -x y z & x y z & -x & -1
\end{array}
\right].
\eeq
Then $B$ has a $2$-circulant representation.
\end{proposition}
\begin{proof}
Let us define permutational  matrices
$P=\left[e_1,e_4,e_2,e_5,e_3,e_6\right],
Q=\left[e_5,e_1,e_3,e_4,e_6,e_2\right]$, and the following unitary
diagonal matrices
$D_1=\mathrm{Diag}\left(1,\sqrt[3]{z},1/\sqrt[3]{z},1/y,\sqrt[3]{z},-1/(x
y \sqrt[3]{z})\right)$ and
$D_2=\mathrm{Diag}\left(1,1/\sqrt[3]{z},1/z^{2/3},1,-x y z^{2/3},y
\sqrt[3]{z}\right)$. Here, $\sqrt[3]{z}$ denotes the principal
cubic root of $z$, and $z^{2/3}$ is the (slightly abusive)
notation of $\left(\sqrt[3]{z}\right)^2$. Now we see that
$D_1PBQD_2$ is $2$-circulant, in particular \beql{31re}
D_1PBQD_2=\left[
\begin{array}{ccc|ccc}
 1 & \frac{1}{\sqrt[3]{z}} & \frac{1}{z^{2/3}} & 1 & -x y z^{2/3} & y
   \sqrt[3]{z} \\
 \frac{1}{z^{2/3}} & 1 & \frac{1}{\sqrt[3]{z}} & y \sqrt[3]{z} & 1 & -x y
   z^{2/3} \\
 \frac{1}{\sqrt[3]{z}} & \frac{1}{z^{2/3}} & 1 & -x y z^{2/3} & y \sqrt[3]{z}
   & 1 \\
   \hline
 1 & \frac{1}{y \sqrt[3]{z}} & -\frac{1}{x y z^{2/3}} & -1 & -z^{2/3} &
   -\sqrt[3]{z} \\
 -\frac{1}{x y z^{2/3}} & 1 & \frac{1}{y \sqrt[3]{z}} & -\sqrt[3]{z} & -1 &
   -z^{2/3} \\
 \frac{1}{y \sqrt[3]{z}} & -\frac{1}{x y z^{2/3}} & 1 & -z^{2/3} &
   -\sqrt[3]{z} & -1
\end{array}
\right].
\eeq
\end{proof}
As the elegant characterization of Beauchamp and Nicoara \cite{BN}
shows, all self-adjoint Hadamard matrices of order $6$ are
equivalent to a matrix described by \eqref{ccc}.
\begin{corollary}
All self-adjoint Hadamards of order $6$ has the $2$-circulant
representation.
\end{corollary}
We close this section with the following remark: matrices $M_6$ and $S_6$ are not members of the family $X_6(\alpha)$. It was explicitly stated in \cite{MS}, that $M_6$ and $\overline{M}_6$ are inequivalent, and hence a local neighborhood around the one-parametric matrix $M_6$ avoids the family $X_6$, which is stable under conjugation. Clearly, as $S_6$ is isolated, it cannot be a member of a continuous family of matrices.
\section{The existence of a two-parameter family of MUB-triplets in $\mathbb{C}^6$}
Recall that a family of mutually unbiased  bases (MUBs)
$\{\mathcal{B}_1,\mathcal{B}_2,\hdots,\mathcal{B}_k\}$ is a
collection of orthonormal bases of $\mathbb{C}^n$ such that
$\left|\left\langle e,f\right\rangle\right|=1/\sqrt{n}$ whenever
$e\in B_i$ and $f\in B_j$ for some $i\neq j$. One can assume that
$\mathcal{B}_1$ is the standard basis, and hence the coordinates
of the vectors of all the remaining bases have modulus
$1/\sqrt{n}$. In particular, the column vectors of the remaining
bases --- up to a constant factor
--- form complex Hadamard matrices. It is well known that at most
$n+1$ MUBs can be constructed, and this upper bound is sharp
whenever $n$ is a prime power. On the other hand, when $n$ is
composite, not much is known about the existence of mutually
unbiased bases. For a quick introduction to MUBs we refer the
reader to \cite{MUB6}, \cite{wer}. In this section the existence
of a two-parameter family of {\it MUB-triplets} of order $6$ is
concluded. The method described here was discovered by Zauner
\cite{Z}, who exhibited a one-parameter family of triplets
earlier.\footnote{The existence of a one-parameter family of
MUB-triplets of order $6$ was also discovered very recently in
\cite{uj} by a method completely different from Zauner's
approach.} Interestingly, the heart of his construction was the
existence of the infinite family of $2$-circulant complex Hadamard
matrices described by formula \eqref{f32}, which he used as a
\emph{seed matrix} for producing MUB triplets. We recall his
machinery and apply it to the two-parameter matrix $X_6(\alpha)$.
First we recall a simple, but extremely useful lemma on the
representation of $2\times 2$ unitaries.
\begin{lemma}[cf.\ Lemma 5.5.\ from \cite{Z}]\label{ZL}
Suppose that $M$ is a $2\times 2$ unitary matrix with entries $a,b,c$ and $d$. Then there exists $u,v,x,y\in \mathbb{T}$, such that
\beql{}
M=\left[\begin{array}{cc}
a & b\\
c & d\\
\end{array}\right]=
\frac{1}{2}\left[\begin{array}{rr}
u+v & y\left(u-v\right) \\
\frac{\left(u-v\right)}{x} & \frac{y\left(u+v\right)}{x}\\
\end{array}\right].
\eeq
$\hfill\square$
\end{lemma}
Before proceeding we need to introduce some notations. Let $T$ be a $2m\times 2m$ block matrix with $m\times m$ blocks as the following:
\beql{}
T=\left[
\begin{array}{cc}
A  & B\\
C & D\\
\end{array}
\right]. \eeq Further let $U,V,X,Y$ are arbitrary unitary diagonal
matrices, and let us define the following matrices with the aid of
the Fourier matrix $F_m$ as \beql{z1z2} Z_1=\frac{1}{\sqrt2}\left[
\begin{array}{rr}
F_m & XF_m\\
F_m & -XF_m\\
\end{array}
\right],
Z_2=\frac{1}{\sqrt2}\left[
\begin{array}{rr}
UF_m & UYF_m\\
VF_m & -VYF_m\\
\end{array}
\right]. \eeq Note that as $F_m$ is unitary, so are $Z_1, Z_2$ and
hence also \beql{zi}
Z_1^{-1}Z_2=\frac{1}{2}\left[\begin{array}{cc}
F_m^{-1}\left(U+V\right)F_m & F_m^{-1}\left(\left(U-V\right)Y\right)F_m\\
F_m^{-1}\left(X^{-1}\left(U-V\right)\right)F_m & F_m^{-1}\left(X^{-1}Y\left(U+V\right)\right)F_m\\
\end{array}
\right]. \eeq In \cite{Z} Zauner characterized $2$-circulant
unitary matrices in the following way. We quote his result with a
sketched proof for completeness.
\begin{proposition}[cf.\ Prop.\ 5.6.\ from \cite{Z}]\label{ZP}
$T$ is a $2$-circulant unitary matrix with blocks $A, B, C, D$ if
and only if there exist $2m\times 2m$ (rescaled) complex Hadamard
matrices $Z_1, Z_2$ as in formula \eqref{z1z2}, such that
$T=Z_1^{-1}Z_2$.
\end{proposition}
\begin{proof}[Proof (Sketch)]
Suppose that $Z_1, Z_2$ are given as above. Clearly $Z_1^{-1}Z_2$
is unitary. Also, for every diagonal matrix $D$ the
matrix $F_m^{-1}DF_m$ is circulant, and hence $T$ is a
$2$-circulant unitary by formula \eqref{zi}.

For the converse, suppose that $T$ is an arbitrary $2$-circulant
unitary matrix. Then one can write $T$ as \beql{zi2}
T=\left[\begin{array}{cc}
F_m^{-1} & 0\\
0 & F_m^{-1}\\
\end{array}\right]\left[\begin{array}{rr}
\tilde{A} & \tilde{B}\\
\tilde{C} & \tilde{D}\\
\end{array}
\right]\left[\begin{array}{cc}
F_m & 0\\
0 & F_m\\
\end{array}\right],
\eeq with diagonal matrices $\tilde A, \tilde B, \tilde C, \tilde
D$. It follows that $T$ is unitary if and only if the matrices
\beql{} S_k=\left[
\begin{array}{rr}
a_{k} & b_{k}\\
c_{k} & d_{k}\\
\end{array}
\right] \eeq are unitary for every  $1\leq k\leq m$. Now use Lemma
\ref{ZL} to represent $S_k$ with unimodular elements $u_k, v_k,
x_k, y_k$, from which one readily defines the unitary diagonal
matrices $U, V, X, Y$, and finally $Z_1$ and $Z_2$ through formula
\eqref{z1z2}. We conclude by observing that in this setting
formulas \eqref{zi} and \eqref{zi2} coincide.
\end{proof}
Proposition  \ref{ZP} describes how to construct a triplet of MUBs
from a given $2$-circulant complex Hadamard matrix $T$. Clearly,
the assumption that $T=Z_1^{-1}Z_2=Z_1^{\ast}Z_2$ is a complex
Hadamard matrix implies that $\{I,Z_1,Z_2\}$ is a collection of
$3$ MUBs of order $2m$. Note, however, that in order to use this
construction one needs to begin with a suitable complex Hadamard
matrix $T$ first, from which the unbiased bases $Z_1$ and $Z_2$
can be constructed. Clearly, the newly discovered matrix
$X_6(\alpha)$ is a perfect seed matrix for Zauner's construction.
In summary, we have proved the following
\begin{theorem}
There exists a two-parameter family of MUB-triplets of order $6$
emerging from the family $X_6(\alpha)$ via Zauner's construction
described in Proposition \ref{ZP}.
\end{theorem}
We conclude our paper by the  following observation: it is
plausible that our new family $X_6(\alpha)$ intersects the Fourier
families only at $\alpha=0$. If one could exhibit similar families
$X_6^{{a,b}}(\alpha)$ for all members $F(a,b)$ of the Fourier
families, that would provide a rigorous proof of the existence of
a 4-parameter family of complex Hadamard matrices of order $6$.
The existence of such a family is strongly indicated by the
numerical results of \cite{skinner}. This could possibly lead to a
full classification of complex Hadamard matrices of order $6$.
\section*{Acknowledgement}
The author thanks Ingemar Bengtsson, M\'at\'e Matolcsi and Karol \.Zyczkowski for various useful remarks concerning this manuscript.

\end{document}